\documentclass[12pt,a4paper]{amsart}
\usepackage{amssymb,amsmath,amsthm}
\usepackage{graphicx}
\usepackage{enumerate}
\usepackage{tikz}

\long\def\onefigure#1#2{
\begin{figure*}[tbp]
\begin{center}
#1
\end{center}
\caption{#2}
\end{figure*}
} 

\newcommand{\lipefig}[2]  
{\onefigure{\mbox{\psfig{file=#1.eps}}}{\label{f:#1} #2} }

\newtheorem{theorem}{Theorem}[section]
\newtheorem{lemma}[theorem]{Lemma}
\newtheorem{corollary}{Corollary}[section]

\newtheorem{claim}{Claim}[section]
\newtheorem{fact}{Fact}[section]

\newcommand{\remove}[1]{}

\newcommand{\de}{\delta}
\newcommand{\al}{\alpha}
\newcommand{\be}{\beta}
\newcommand{\ga}{\gamma}

\newcommand{\si}{\sigma}
\newcommand{\eps}{\varepsilon}

\newcommand{\rr}{\mathbb{R}^d}
\newcommand{\R}{\mathbb{R}}

\newcommand{\conv}{\mathrm{conv}}

\newcommand{\lin}{\textrm{lin\;}}
\newcommand{\dist}{\textrm{dist}}

\numberwithin{equation}{section}

\begin{document}

\title{Universal sequences of lines in $\R^d$}
\author{Imre B\'ar\'any, Gil Kalai, and Attila P\'or}


\maketitle

\begin{abstract} 
One of the most important and useful examples in discrete geometry is a finite sequence of points on the moment curve $\gamma(t)=(t,t^2,t^3,\dots ,t^d)$ or, more generally, on a {\it strictly monotone curve} in $\mathbb R^d$. These sequences  as well as the ambient curve itself can be described in terms of {\it universality properties} and we will study the question:  "What is a universal sequence of oriented and unoriented lines in $d$-space''

We give partial answers to this question, and to the analogous one for $k$-flats. Given a large integer $n$, it turns out that, like the case of points the number of universal configurations is bounded by a function of $d$, but unlike the case for points, there are a large number of distinct universal finite sequences of lines. We show that their number is at least $2^{d-1}-2$ and at most $(d-1)!$. However, like for points, in all dimensions except $d=4$, there is essentially a unique
{\em continuous} example of a universal family of lines. The case $d=4$ is left as an open question.  
\end{abstract}

\section{Introduction}\label{sec:introd}

The vertices of the cyclic polytope in $\mathbb R^d$ is one of the most important examples in discrete geometry. To define them we consider $n$ points on the moment curve $\gamma(t)=(t,t^2,t^3, \dots ,t^d)$ in $\mathbb R^d$, and we can consider any {\it strictly monotone curve,} namely a curve with the property that every $d+1$ points on it are in general position, in place of the moment curve. 
We study high dimensional analogs where "points'' are replaced with lines or even $k$-flats, $k\ge 1$. Our stating point is the universality property of $n$ points in cyclic position and we explore both finite and continuous sets of lines (and $k$-flats) with similar universality property.

Assume $a_1,\ldots, a_N$ is a sequence of points in $\rr$ in general position meaning that no $d+1$
of these $N\ge d+1$ points lie on a hyperplane.
In other words, for every subsequence $1\le i_1 < \ldots < i_{d+1}\le N$ the determinant of the $(d+1)\times (d+1)$ matrix
\begin{equation}
\begin{pmatrix}
a_{i_1} & a_{i_2}& \cdots  & a_{i_{d+1}}\\
1          & 1        &1  \cdots  & 1
\end{pmatrix}
\end{equation}
is different from zero. Write $\det(a,i_1,\ldots,i_{d+1})$ for this determinant. Observe that the general position condition simply means that
the $d+1$ tuples of the sequence are outside the zero set of $m={N \choose d+1}$ polynomials $p_1,\ldots,p_m$.
These polynomials split the set of sequences
of $N$ points into finitely many cells, where a cell is the set of point sequences
of length $N$ such that, for every $i \in [m]$, the sign of $p_i$ is constant, $+1$ or $-1$, on the sequence.

The $m$ signs of these polynomials are referred to as the {\it order-type} of the sequence $a_1,a_2,\dots,a_N$. Moreover, it is known that these  polynomials split the sets of sequences of $N$ points into finitely many connected components as well; their number is finite according to a famous theorem of Oleinik--Petrovskii~\cite{Oleinik49}, Milnor~\cite{Milnor64}, and Thom~\cite{Thom65}.
As a matter of fact, the upper bounds given by these authors imply upper bounds on
the number of order-types, much below the immediate bound $2^m$.
For more details see Goodman and Pollack \cite {GooPol86} and Alon \cite {Alo86}.
On the other hand, it is known that there is no upper bound on the
number of connected components described by a single order type, and, in fact, a single order type can be as topologically complicated
as essentially any algebraic variety, see Mn\"ev \cite{Mnev88}  and Richter-Gebert \cite{Rich-Geb96}.

A useful and famous (and probably folklore) result says that for every $d$ and $n\ge d+1$ there is an $N=N(d,n)$ such that the following holds.
Every sequence $a_1,\ldots, a_N$ of points in $\rr$ in general position contains a {\sl homogeneous}
subsequence $b_1,\ldots,b_n$ of length $n$; homogeneous meaning that
the determinants $\det(b,i_1,\ldots,i_{d+1})$ have the same sign for all
sequences $1\le i_1 < \ldots < i_{d+1}\le n$. Subsequence means, as usual,
that $b_j=a_{k_j}$ for all $j\in [n]$ where $1\le k_1 < \ldots < k_n\le N$.
Here $[n]$ denotes the set $\{1,2,\ldots,n\}$. This result actually follows
from Ramsey's theorem~\cite{Ramsey29}: the $d+1$ tuples of our sequence are coloured by
$+1$ or by $-1$, so there is a large subsequence all of whose $d+1$ tuples
are of the same colour. More precisely, for every $n$ there is an $N$ such that
every $d+1$ tuple of an $n$ element subsequence of $a_1,\ldots,a_N$ carries the same colour.
When $d=2$ this method gives a proof of the famous Erd\H os--Szekeres theorem,
see \cite{ESz35} and \cite{ESz61}, but with weaker bounds on $N$. We remark that in this paper
an $r$ tuple of points (or objects) always means an ordered $r$ tuple, that is, a sequence of points (or objects) of length $r$.

Assume $0<t_1<\ldots <t_n$ and consider the sequence of points $\ga(t_1),\ldots,\ga(t_n)$ from
the moment curve $\ga(t)=(t,t^2,\ldots,t^d), t \in \R$.
Observe that the determinant $\det(\ga(t),i_1,\ldots,i_{d+1})$ is positive for every sequence
$1\le i_1<\ldots <i_{d+1}\le n$. Set $\ga^*(t)=(t,t^2,\ldots, t^{d-1},-t^d)$.
Analogously for the sequence $\ga^*(t_1),\ldots,\ga^*(t_n)$ the corresponding determinants are all negative.

The property that all determinants have the same sign is {\sl universal}: for every $n$,
every long enough sequence of points in $\rr$ has a subsequence of length $n$
with this property. The two examples with $\ga$ and $\ga^*$ show that there are exactly
two kinds of {\sl universal sequences}. In the language of polynomials and cells,
universality says that our subsequence lies in a cell $C$, that is,
there is a subsequence $b_1,\ldots,b_n$ such that $C$ is on the positive side of each polynomial
 $\det(b,i_1,\ldots,i_{d+1})$, or on the negative side of each such polynomial.

Here comes another example. Again, let $a_1,\ldots, a_N$ be a sequence of points in $\rr$ in general position. Consider a subsequence
$b_1,\ldots,b_{d+2}$ and set $B=\{b_1,\ldots,b_{d+2}\}$. Note $B$ has a unique Radon partition \cite{Radon21}, that is $[d+2]=X \cup Y$ with $X,Y\ne \emptyset$ and
$X,Y$ are disjoint and $\conv \{b_i:i\in X\} \cap \conv \{b_j: j\in Y\} \ne \emptyset$, for concreteness we assume $1\in X$. So each subsequence
$b_1,\ldots,b_{d+2}$ defines a subset $X\subset [d+2]$ (with $1\in X \ne [d+2]$). Ramsey's theorem applies again and implies that, for $N$ large enough,
there is a subsequence $b_1,\ldots,b_n$ of the $a_i$ such that $X$ is the same subset of $[d+2]$ for every subsequence $c_1,\ldots,c_{d+2}$ of the $b_i$ sequence.
What is the universal type or what are the universal types for this {\sl Radon property}? It turns out that the answer is simple: $X$ and $Y$ are interlacing
subsets of $[d+2]$, that is, $X=\{1,3,5,\ldots\}$ and $Y=\{2,4,\ldots\}$. The proof is left to the interested reader. One can check that in the above examples
with $\ga(t)$ and $\ga^*(t)$ the universal sequences indeed have interlacing Radon partitions.

We conclude the introduction by describing the structure of the paper and with it our main results. In Section \ref {sec:univer} we outline the definition of universal sequences of lines and $k$-flats. In Section \ref {sec:lines} we consider oriented lines and show that every universal sequence of oriented lines in $\mathbb R^d$ is described by a permutation on $\{1,2,\dots, d-1\}$. This implies that for every fixed dimension $d$ there is only a finite number of universal sequences of lines. A similar treatment of universal unoriented lines is given in Section \ref {sec:ulines}, and for $k$-flats when $d \equiv 1 (\mod k)$, in Section \ref {sec:kflats}. Following some necessary results on matrices with rapidly increasing entries (Section \ref {sec:stretch}; with one proof left to Section \ref {sec:RI}), we prove in Section \ref {sec:prooftyp} that the number of universal sequences of oriented lines is at least $2^{d-1}-2$. This is a place where our theory differs from the classical theory for points; we demonstrate an exponential number (in the dimension $d$) of distinct sequences of universal lines while for points, universal sequences are essentially unique. We regain uniqueness by moving to continuous sequences of lines:
in Section \ref {sec:cont} we define continuous universal sequences of lines and in Section \ref {sec:unique} prove a uniqueness theorem for them for all dimensions except $d=4$.  
 
 We remark that Sturmfels proved \cite {Stu87} that every universal sequence of points in $\mathbb R^d$ can be extended to a strictly monotone curve. Our results show that for lines this is no longer the case.  

\section{Definition of universality}\label{sec:univer}

What is a universal sequence of lines, and of oriented lines in $\rr$? How many types of them are there? This is the main topic in this article. We have some partial results for the case of lines and also for the same question with $k$-flats in $\rr$. The definition of universality (or rather its metadefinition) requires four conditions or steps.

(i) First we define when an ordered $r$ tuple of $k$-flats is in {\sl general position}. This usually means that they do not lie in the zero set of a finite number of well-defined polynomials that determine cells in the space of $r$ tuples of $k$-flats. Here we also require that there be at least two such cells.

(ii) Next we define a {\sl property} of an ordered $s$ tuple $A_1,\ldots,A_s$ of $k$-flats whose $r$ tuples are in general position, of course $s\ge r$.
This property is a function $F$ on $s$ tuples of $k$-flats that take values in a finite set $M$. Elements of $M$ will be called {\sl types}. In our first example this finite set $M$ is $\{1,-1\}$ and in the second $M$ is the family of all $X \subset [d+2]$ with $1 \in X \ne [d+2]$.

(iii) Assume $B_1,\ldots, B_n$ is a sequence of $k$-flats (oriented or unoriented) in $\rr$ all of whose $r$ tuples are in general position. $F$ maps subsequences of length $s\ge r$ of this sequence to elements of $M$. The sequence $B_1,\ldots, B_n$ is called {\sl homogeneous} relative to {\sl property} $F$ if $F$ maps all of its $s$ tuples to the same element (or type) $m \in M$. A type $m \in M$ is {\sl universal} if, for every $n\ge s$ there is a homogeneous sequence $B_1,\ldots, B_n$ whose $s$ tuples are all mapped to $m$.

(iv) Next comes Ramsey:  For every $n\ge s$ there is $N$ such that the following holds. Assume $A_1,\ldots, A_N$ is a sequence of $k$-flats in $\rr$ all of whose $r$ tuples are in general position. As $M$ is finite, Ramsey's theorem implies the existence of a homogeneous subsequence $B_1,\ldots, B_n$ of the sequence  $A_1,\ldots, A_N$. Of course $N$ depends on $n,r,s$ and the map $F$ as well. The finiteness of $M$ implies further that there is a universal element $m$ in $M$. The corresponding sequences $B_1,\ldots,B_n$ are called {\sl universal} of type $m$ (relative to {\sl property} $F$).

\medskip
The question is what types $m \in M$ are universal, how many of them are there, and what the universal sequences look like.

\section{Oriented lines}\label{sec:lines}

We are going to use some notation from exterior algebra, for instance $u_1 \wedge \ldots \wedge u_d$
is the determinant of the matrix whose columns are the vectors $u_i \in \rr, i\in [d]$, and $u_1 \wedge \ldots \wedge u_{d-1}$ is a vector in $\rr$, the wedge product
of the $u_i$s.

An oriented line $L$ in $\rr$ $d\ge 2$ is given by a pair $(a,v)$ with $a,v \in \rr$ and $v \ne 0$ and $L=\{a+tv:t \in \R\}$. In fact it is the equivalence class
of such pairs where $(a,v)$ and $(\bar{a},\bar{v})$ are equivalent (or represent the same line) if $\bar{a}=a+\al v$ for some $\al \in \R$ and $\bar{v}=\be v$ with $\be>0$.
An ordered $d-1$ tuple of lines $L_1,\ldots,L_{d-1}$, with $L_i$ represented by $(a_i,v_i)$ is in {\sl general position} if the numbers
$h_i:=a_i\wedge v_1 \wedge \ldots \wedge v_{d-1}$, $i\in [d-1]$ are all distinct. General position then means that the system $(a_i,v_i), i \in [d-1]$ avoids
the zero set of the polynomials $(a_i-a_j)\wedge v_1 \wedge \ldots \wedge v_{d-1}$, for distinct $i,j \in [d-1]$.
In particular, $u=v_1 \wedge \ldots \wedge v_{d-1}$ is a non-zero vector in $\rr$. We mention that $h_i=u\cdot v_i$, scalar product.

The $d-1$ real numbers $h_i$ come in increasing order as $h_{j_1}<h_{j_2}< \ldots <h_{j_{d-1}}$. They define a permutation $\si$ of $[d-1]$ via $\si(i)=j_i$.
As is easy to check this permutation $\si$ does not depend on the choice of the pair $(a_i,v_i)$ representing $L_i$: $\si$ depends only on the $d-1$ tuple $L_1,\ldots, L_{d-1}$. Observe that every permutation of $[d-1]$ can occur. Here we need $d\ge 3$ as for $d=2$ there is only one $h$.

Another way to see the permutation $\si$ is to consider the hyperplane $H=\frac 1{d-1}(L_1+\ldots +L_{n-1})$ which is the Minkowski average of the lines.
Because of the general position assumption $H$ is indeed a hyperplane, its outer normal is $u=v_1 \wedge \ldots \wedge v_{d-1}\ne 0$. A suitably
translated copy, say $H_i$ of $H$ contains $L_i$. The hyperplanes $H_1,\ldots,H_{d-1}$ intersect the line whose direction is $u$ in order $\si$
in distinct points.

This time $F$ maps $L_1,\ldots, L_{d-1}$ to $\si$ and the set of permutations is finite, of size $(d-1)!$. Ramsey's theorem applies and gives the following.

\begin{theorem}\label{th:olines} For integers $n \ge d \ge 3$ there is a number $N$ such that every sequence $L_1,\ldots,L_N$ of oriented lines in $\rr$ whose $d-1$ tuples are in general position contains a homogeneous subsequence $K_1,\ldots,K_n$ of type $\si$ for some permutation $\si$ of $[d-1]$. In other words $F$ maps every $d-1$ tuple of $K_1,\ldots,K_n$ to the same type $\si$.\qed
\end{theorem}

\begin{corollary}\label{cor:olines} For oriented lines in $\rr$ there are universal types, that is permutations  $\si$ of $[d-1]$. There are at most $(d-1)!$ of them.\qed
\end{corollary}

\begin{figure}[h!]
\centering
\includegraphics[scale=0.7]{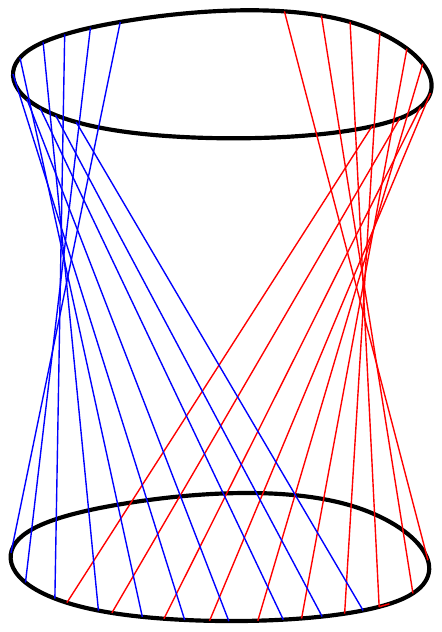}
\caption{Two sets of universal oriented lines in $\R^3$.}
\label{fig:hyperb}
\end{figure}

This implies in particular, that in $\R^3$ there are at most two types. In fact there are exactly two types in $\R^3$. This can be seen from the example of the one sheeted hyperboloid whose equation is $x^2+y^2=z^2+1$. This hyperboloid contains two sets of lines, see Figure~\ref{fig:hyperb}. Given a line $L$ on the hyperboloid, it is associated with the pair $(a,v)$. Orient $L$ by requiring that the $z$ component of $v$ is positive. The two sets of lines are shown in Figure~\ref{fig:hyperb}.

The question is how many of the possible $(d-1)!$ types are universal. A partial answer is the content of the next result.

\begin{theorem}\label{th:types} There are at least $2^{d-1}-2$ different universal permutations of $[d-1]$.
\end{theorem}

The proof is in Section~\ref{sec:prooftyp}. It uses rapidly increasing sequences and RI matrices that are related to the stretched grid and to the stretched diagonal, that come from a paper by Bukh, Nivasch,  Matou\v{s}ek~\cite{BukhMN11}, and are also connected to a construction of P\'or~\cite{Por18}. The necessary background is given in Section~\ref{sec:stretch}.

We remark here that the red lines in Figure~\ref{fig:hyperb} form a {\sl continuous and universal family} of lines. More precisely, define $a(t)=(\cos t, \sin t, 0)$ and $v(t)=(-\sin t,\cos t, 1)$. Let $L(t)$ be the oriented line given by the pair $(a(t),v(t)), t \in [0,\pi)$. When $0 < t_1<t_2 <\ldots <t_n <\pi$ the sequence of lines $L(t_1),\ldots,L(t_n)$ is homogeneous of type identity, as one can check directly. That's why we call the family of lines $L(t), t \in [0,\pi)$ a continuous and universal family of lines.  Section~\ref{sec:cont} gives the proper definition, and examples of such families in every dimension. In Section~\ref{sec:unique} we prove that the type of such a family is either the identity or its reverse.

\section{Unoriented lines}\label{sec:ulines}

\begin{figure}[h!]
\centering
\includegraphics[scale=0.8]{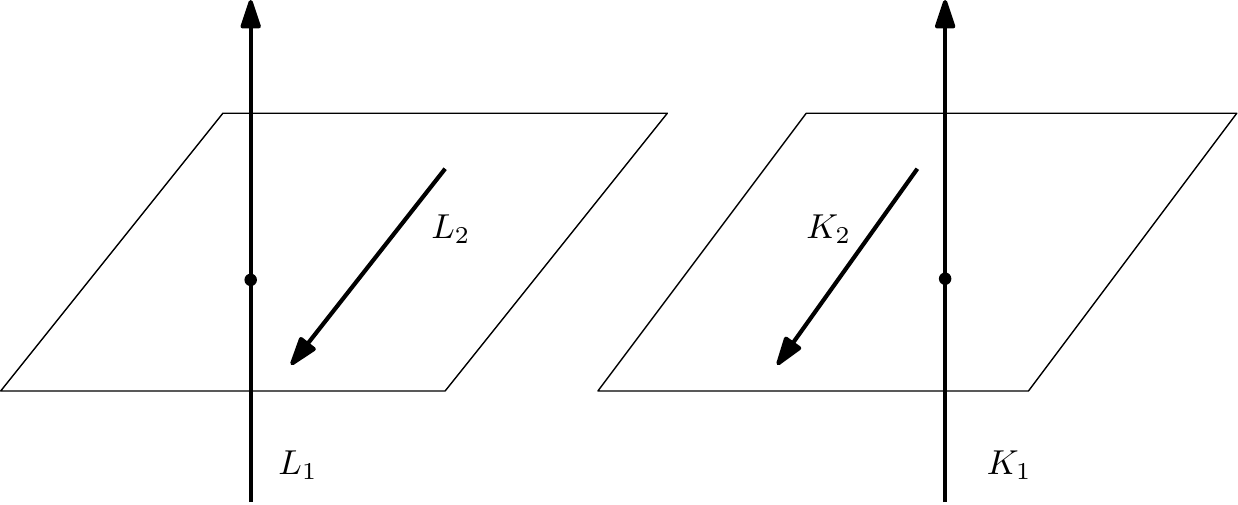}
\caption{Non-homotopic pairs of oriented lines in $\R^3$.}
\label{fig:no-homo}
\end{figure}

The setting with unoriented lines (or simply lines) in $\rr$ is similar to the oriented ones. But for instance in $\R^3$ there are pairs $(L_1,L_2)$ and $(K_1,K_2)$ of oriented lines in general position that belong to distinct cells (defined in this case by a single polynomial), see Figure~\ref{fig:no-homo}, they are of distinct types. This does not hold for lines: these pairs (when unoriented) belong to the same cell. In fact, they can be carried to each other by a homotopy through general position pairs. This is not the case for triples
of lines. That is, there are triples $(L_1,L_2,L_3)$ and $(K_1,K_2,K_3)$ of lines in $\R^3$ (in general position) such that belong to distinct cells, as we shall see soon.

A line $L$ in $\rr$ is given by a pair $(a,v)$ with $a,v \in \rr$  and $v \ne 0$ and $L=\{a+tv:t \in \R\}$. In fact it is the equivalence class
of such pairs where $(a,v)$ and $(\bar{a},\bar{v})$ are equivalent (or represent the same line) if $\bar{a}=a+\al v$ for some $\al \in \R$ and $\bar{v}=\be v$ with $\be\ne 0$. Given an ordered $d$ tuple of lines $L_1,\ldots,L_d$, with $L_i$ represented by $(a_i,v_i)$, define $u_j=v_1 \wedge \ldots \wedge v_d$ where $v_j$ is missing from the wedge product. The $d$ tuple $L_1,\ldots,L_d$ is in {\sl general position} if, for every $j \in [d]$, the numbers
$h_{j,i}:=a_i\wedge v_j$, $i\in [d]\setminus \{j\}$ are all distinct. General position again means that the system $(a_i,v_i), i \in [d]$ avoids
the zero set of certain polynomials. Note that we assume here $d>2$.

For every $j \in [d]$ the numbers $h_{j,i}, i \in [d]\setminus \{j\}$ define a permutation $\si_j$ of $[d]\setminus \{j\}$, namely,
$h_{j,\si_j(1)}<h_{j,\si_j(2)}< \ldots < h_{j,\si_j(d)}$ where again $h_{\si_j(j)}$ is not defined and is missing from the list.
Note that $L_i$ is represented by both $(a_i,v_i)$ and $(a_i,-v_i)$. Consequently $\si_j=(\si_j(1),\ldots,\si_j(d))$ and $\si_j^*=(\si_j(d),\ldots,\si_j(1))$
represent the same ordering, of course $\si_j(j)$ is missing again.

So the map $F$ for universality associates with $L_1,\ldots,L_d$ $d$ pairs of permutations $\{\si_1,\si_1^*\},\ldots,\{\si_d,\si_d^*\}$. The values of $F$ are from a finite set $M$, of size $((d-1)!/2)^d$, so Ramsey's theorem works again:

\begin{theorem}\label{th:ulines} For integers $n,d$ with $n \ge d \ge 3$ there is a number $N$ such that every sequence $L_1,\ldots,L_N$ of lines in $\rr$ whose $d$ tuples are in general position contains a homogeneous subsequence $K_1,\ldots,K_n$. In other words, $F$ maps every $d$ tuple $K_{i_1},\ldots,K_{i_d}$ to the same type $\{\si_1, \si_1^*\},\ldots, \{\si_d,\si_d^*\}$.\qed
\end{theorem}

A direct corollary is that there are universal types for unoriented lines in $\rr$ and that there are universal sequences of unoriented lines of length $n$ for every $n \ge d$. Their number is at most $((d-1)!/2)^d$. We are going to reduce this number to $(d-1)!/2$. Observe first that $\si_j$ is a permutation of $[d]\setminus \{j\}$ which is a linearly ordered set of $d-1$ elements. So we can consider $\si_j$ a permutation of $[d-1]$, and $\si_j^*$ is its reverse permutation.

\begin{lemma}\label{l:reduce} Assume $L_1,\ldots,L_n$ is a homogeneous sequence of lines in $\rr$, $n>d\ge 3$ with permutation pairs $\{\si_j,\si_j^*\}$ of $[d-1]$ for every $j\in [d]$. Then $\{\si_1,\si_1^*\}=\{\si_2,\si_2^*\}=\ldots=\{\si_d,\si_d^*\}$.
\end{lemma}

{\bf Proof.} It suffices to consider $n=d+1$. We show that $\si_1=\si_d$. For the ordered $d$ tuple $L_1,\ldots,L_d$ $u_1=v_2\wedge \ldots \wedge v_d$, and the permutation $\si_1$ (of $[d-1]$) is determined by the increasing rearrangement of the numbers $a_2\wedge u_1, a_3\wedge u_1,\ldots ,a_d\wedge u_1$.
We check the permutation $\si_d$ for the ordered $d$ tuple $L_2,\ldots,L_{d+1}$. The corresponding $u$ vector is exactly the previous $u_1=v_2\wedge \ldots \wedge v_d$  because the last vector $v_{d+1}$ has to be deleted. So $\si_d$ is given by the increasing rearrangement of the same numbers $a_2\wedge u_1, a_3\wedge u_1,\ldots ,a_d\wedge u_1$. Then $\si_1=\si_d$ indeed. The proof of the other cases $\si_1=\si_j$, $j<d$ is identical.\qed

\begin{theorem} The number of universal types for unoriented lines in $\rr$ is at least $2^{d-2}-1.$
\end{theorem}

{\bf Remark.} Assume that $L_1,\ldots,L_n$ ($n \ge d$) is a universal sequence of oriented lines of type $\si$, a permutation of $[d-1]$. Forgetting their orientation, the same sequence of (unoriented) lines becomes a universal sequence of unoriented lines. Indeed, every $d-1$ tuple of the lines is ordered by $\si$ or $\si^*$ as one can check directly.

\section{Universal $k$-flats when $d \equiv 1$ mod $k$}\label{sec:kflats}

Next we consider oriented $k$-flats in $\rr$ under the condition that $d \equiv 1$ mod $k$, or, in different form $d=rk+1$ with $r\ge 2$ and integer. In this case the method we used for lines works so we only give a sketch. An oriented $k$-flat $A$ is given by a pair $(a,B)$ where $a\in \rr$ and $B$ is an ordered set of $k$ linearly independent vectors: $B=(v_1,\ldots,v_k)$. The linear span of $B$, $\lin B$, is defined as
\[
\lin B=\left\{\sum_1^k \al_iv_i: \al_i \in \R \mbox{ for all } i\in [k]\right\}.
\]
Then $A=a+\lin B$. In fact $A$ is given by an equivalence class of such pairs where $(a,B)$ and  $(\bar{a},\bar{B})$ are equivalent if $\bar{a}=a+\sum_1^k \al_iv_i$ with some real numbers $\al_i$, and if there is a linear transformation $T: \lin B \to \lin B$ with positive determinant that carries the basis of $B$ to the basis of $\bar{B}$, that is $Tv_i=\bar{v_i}$ for $i \in [k],$ where $\bar{B}=(\bar{v_1},\ldots,\bar{v_k}).$
Define  $u(B)= v_1\wedge \ldots \wedge v_k.$

Given an $r$ tuple $A_1,\ldots,A_r$ of $k$ flats with $A_i$ represented by $(a_i,B_i)$ set
\[
u=u(B_1) \wedge \ldots \wedge u(B_r).
\]
The condition $d \equiv 1$ mod $k$ implies that $u$ is a vector in $\rr$. The $r$ tuple $A_1,\ldots,A_r$ is in {\sl general position} if the numbers $h_i:=a_i \wedge u$, $i\in [r]$ are all distinct. General position then means that the system $(a_i,B_i), i \in [r]$ avoids the zero set of certain polynomials: $(a_i-a_j)\wedge u$, for distinct $i,j \in [r]$, in particular $u\ne 0$.
The increasing rearrangement of the numbers $h_1,\ldots,h_r$ defines a permutation $\pi$ of $[r]$. We observe that $\pi$ depends only on the $r$ tuple $A_1,\ldots,A_r$ and not on the representation by $(a_i,B_i)$ of the $k$-flats $A_i$.

Same way as before $F$ maps the $r$ tuple $A_1,\ldots,A_r$ to the permutation $\pi$, and the universality scheme gives the following result.

\begin{theorem}\label{th:flats} Assume $d=kr+1$ where $d,k,r$ are integers, $r\ge 2, k\ge 1$. There is a universal permutation $\pi$ of $[r]$. This means that for every integer $n>r$ there is a number $N$ such that every sequence $A_1,\ldots,A_N$ of oriented $k$-flats in $\rr$ whose $r$ tuples are in general position contains a subsequence $D_1,\ldots,D_n$ such that $F$ maps every $r$ tuple $D_{i_1},\ldots,D_{i_r}$ to the permutation $\pi$.\qed
\end{theorem}

So under the above conditions there are universal sequences of $k$-flats. Evidently, there are at most $r!$ different types of them. A straightforward
modification of the proof of Theorem~\ref{th:types} shows that their number is at least $2^r-2$.

\section{Rapidly increasing entries}\label{sec:stretch}

Let $M$ be a $D \times m$ matrix with entry $a(i,j)$ in row $i$ and column $j$, we assume $D$ is fixed and $m$ is large, much larger than $D$. The entries in $M$ are {\sl rapidly increasing} if every $a(i,j)\ge 1$ is an integer and, for fixed $i \in [D]$, $a(i,j+1)$ is much larger than $a(i,j)$ for all $j\in [m-1]$,
and further, $a(i+1,1)$ is larger than $a(i,m)$. Such a matrix is called an RI matrix. Similar matrices with various purposes were constructed by Bukh, Nivasch, and Matou\v{s}ek~\cite{BukhMN11} and by P\'or~\cite{Por18}.

The main feature of an RI matrix is that if $M^*$ is a $k\times k$ submatrix of $M$ (here $k\in [D]$ and $k \ge 2$), then $\det M^*$ is essentially equal to the product of the entries on the main diagonal of $M^*$. The meaning of ``essentially equal" is made precise the following way.
Given a small $\eps>0$, the ``much larger" in the definition of an RI matrix can be chosen so large that, with $P$ denoting the product of the entries
on the main diagonal of $M^*$, we have $|\det M^*-P|<\eps P.$ With a slight but very convenient abuse of notation we will write this as
\begin{equation}\label{eq:RI}
\det M^*=(1\pm \eps)P.
\end{equation}

The ``much larger" condition in the definition of rapidly increasing is in fact $Dm-1$ conditions.
It follows from the results of \cite{BukhMN11} and \cite{Por18} that they can be chosen so that equation~(\ref{eq:RI}) holds.
We will need one further requirement, namely

\begin{eqnarray}\label{eq:extraC}
&a(i,j-1)a(i,j)<\eps a(i,j+1)\\
 &\mbox{ for all }i \in [D]\mbox{ and }j\in \{2,\ldots,m-1\}.\nonumber
\end{eqnarray}
\begin{theorem}\label{th:RI} For integers $D,m$ with $m>D\ge 1$ and for every $\eps>0$ there is an $D \times m$ RI matrix $M$
satisfying (\ref{eq:RI}) and (\ref{eq:extraC}).
\end{theorem}

The proof is postponed to Section~\ref{sec:RI}. The extra condition (\ref{eq:extraC}) implies that $a(i,j)<\eps a(i,j-1)a(i,j+1)$
for all $i \in [D]$ and $j\in \{2,\ldots,m-1\}.$ This follows from the inequality
\[
a(i,j)< \frac {\eps a(i,j+1)}{a(i,j-1)}\le  \eps a(i,j+1) a(i,j-1)
\]
as $a(i,j-1)\ge 1$. We will use the following consequences of the extra condition.
For all $i \in [D]$
\begin{eqnarray}\label{eq:extra}
a(i,j) &<& \eps a(i,j-1)a(i,J) \mbox{ when }J>j,  \mbox{ and }\\
a(i,J)a(i,j)&<&\eps a(i,j+1)\mbox{ when }J<j.\nonumber
\end{eqnarray}

\section{Proof of Theorem~\ref{th:types}}\label{sec:prooftyp}

We are going to give examples of universal permutations other than the identity and the reverse identity. These examples are number sequences generated by writing the numbers $1, \ldots , d-1$ one by one, starting by $1$ and appending the next number either on the left or the right of the current sequence, for instance
$$9, 8, 5, 4, 1, 2,3,6,7, 10$$
or the same but starting from $d-1$ and going downwards. We call these permutations {\sl two-sided stacked}. It is easy to see that there are at most $2^{d-1}-2$ of them.
Indeed, starting with $1$ we have $2^{d-2}$ choices to go left or right, and the same number when starting with $d-1$. But the sequence $1,2,\ldots,d-1$ and its reverse are counted twice.

If $v$ is a $d$ dimensional vector let $[v]_i$ denote its $i$th component. Let $M$ be $2d \times m$ RI matrix satisfying conditions (\ref{eq:extraC}). The parameter $\eps>0$ will be specified later. A typical column $C$ of $M$ is a $2d$ dimensional vector, to be denoted by $C=(v,B)$  where $v,B$ are $d$-dimensional vectors and the ordering within $C$ is as follows. The components of $v$ and $B$ alternate: $[v]_1 < [B]_1 < [v]_2 < \ldots < [B]_d$ where $<$ means in fact``much larger" because of the RI condition.

Define the vector $b\in \rr$ by
\[
[b]_i = \prod_{j=d+1-i}^d [B]_j  \mbox{ so } \frac{[b]_i}{[b]_{i-1}} = [B]_{d+1-i}.
\]

We will choose a sequence $\ga_1,\ldots,\ga_m$ of positive reals, to be specified later, that grow faster than anything else so far. When $v$ is the vector in the column $C=(v,B)$ of $M$ we write $(v)_{-k}$ for the corresponding vector $k$ position before the column of $v$. Let $\delta = (\delta_1, \ldots, \delta_d)$ be a $\pm 1$ vector, that is each $\delta_i$ is either $1$ or $-1$. For a fixed column $C=(v,B)$ of $M$ define $a = \gamma(\sum_{j=1}^d \delta_j [b]_j (v)_{-j})$. This defines the oriented line
corresponding to the pair $(a,v).$ We will only consider a subset of these lines.

Assume $m=n(d+1)$, $n\ge d$ are integers. The sequence in our example consists of $n$ oriented lines, each corresponding to a column
of the form $C_{x(d+1)}$, $x \in [n]$. Call such a column {\sl special}. Select $d-1$ special columns and let $v_1,\ldots,v_{d-1}$
be the $v$ vectors (in increasing order) of the selected columns. We will only consider the following $(d-1)(d+1)$ vectors that come in this order (but not necessary consecutively) in $M$:
\[
(v_1)_{-d}, \ldots, (v_1)_{-1}, v_1, (v_2)_{-d}, \ldots, (v_2)_{-1},v_2, \ldots ,v_{d-1}.
\]
The line $L_i$ corresponds to the pair $(a_i,v_i)$, $i \in[d-1].$

We have to estimate how large $\det (a_i,v_1, \ldots, v_{d-1})$ is. Define
\[
P=[v_1]_1 \cdot \ldots \cdot [v_{i-1}]_{i-1} \cdot [v_i]_{i+1} \cdot \ldots \cdot [v_{d-1}]_{d}.
\]
We are going to show first that, with the convenient notation introduced in (\ref{eq:RI}),
\begin{eqnarray}\label{eq:det}
&&\det(a_i,v_1,\ldots,v_{d-1})=\\
&&\;\;\;(1\pm(d+1)\eps) \ga_iP (-1)^{i-1}\delta_{d+1-i}[b_i]_{d+1-i} [(v_i)_{-(d+1-i)}]_i.\nonumber
\end{eqnarray}

Using the definition of the vector $a_i$ we see that
\begin{eqnarray*}
&\det (a_i,v_1, \ldots, v_{d-1}) = \gamma_i \sum_{j=1}^d \delta_j [b_i]_j \det( (v_i)_{-j} , v_1, \ldots, v_{d-1})\\
     &=\gamma_i \sum_{j=1}^d \delta_j [b_i]_j(-1)^{i-1}\det(v_1, \ldots,v_{i-1},(v_i)_{-j},v_i,\ldots, v_{d-1}).
\end{eqnarray*}
The properties of the RI matrix imply that the last determinant is essentially equal to the product
of the entries on the main diagonal of the corresponding matrix. This product equals
$$
[v_1]_1 \cdot \ldots \cdot [v_{i-1}]_{i-1} \cdot [(v_i)_{-j}]_i \cdot [v_i]_{i+1} \cdot \ldots \cdot [v_{d-1}]_{d}=P\cdot  [(v_i)_{-j}]_i
$$
so the product $P$ is a common factor here, implying that
\[
\det( a_i , v_1, \ldots, v_{d-1})=\ga_iP(-1)^{i-1}  \sum_{j=1}^d\delta_j [b_i]_j [(v_i)_{-j}]_i.
\]

Consider $i$ fixed and set $T_j= [b_i]_j [(v_i)_{-j}]_i$, so the last sum is $\sum_{j=1}^d\delta_jT_j$.
Which is the dominant term here? We claim that it is $T_{d+1-i}.$ We are going to show this by proving that for $j\le d+1-i$
\begin{equation}\label{eq:increas}
T_{j-1}< \eps T_j,
\end{equation}
and for  $j> d+1-i$
\begin{equation}\label{eq:decreas}
T_j<\eps T_{j-1}
\end{equation}

Indeed for $j\le d+1-i$ we have
$$
\frac{[b_i]_j}{[b_i]_{j-1}} = [B_i]_{d+1-j} \ge [v_i]_i.
$$
Because of (\ref{eq:extra}) we have $[(v_i)_{-(j-1)}]_i < \eps [v_i]_i {[(v_i)_{-j}]_i}$ implying that
$[b_i]_{j-1} [(v_i)_{-(j-1)}]_i < \eps [b_i]_j [(v_i)_{-j}]_i,$ which is exactly (\ref{eq:increas}).
Similarly, for $j> d+1-i$ we have
$$
\frac{[b_i]_j}{[b_i]_{j-1}} = [B_i]_{d+1-j} \le [B_i]_{i-1} \le  [v_{i-1}]_i.
$$
In view of (\ref{eq:extra}) ${[(v_i)_{-(j-1)}]_i} [v_{i-1}]_i \le \eps [(v_i)_{-j}]_i$ implying that
$[b_i]_j [(v_i)_{-j}]_i \le \eps [b_i]_{j-1} [(v_i)_{-{j-1}}]_i.$ This is again the same as  (\ref{eq:decreas}).

So the dominant term in the sum $\sum_{j=1}^d \delta_j [b_i]_j [(v_i)_{-j}]_i$ is the one $j=d+1-i$,
and we have $T_j\le \eps T_{d+1-i}$ for all $j\ne d+1-i$, each $T_j \ge 1$ of course. It follows that
\[
\sum_{j=1}^d \delta_j [b_i]_j [(v_i)_{-j}]_i=(1\pm(d-1)\eps)\de_{d+1-i}T_{d+1-i}
\]
extending the notation of (\ref{eq:RI}). We choose now $\eps<\frac 1{K(d+1)}$ with $K$ large, $K=100$ or $1000$, say. It is easy to check that equation~(\ref{eq:det}) holds true. More importantly, with this choice of  $\eps$ the factor $(1\pm(d+1)\eps)$ is between $(1-\frac 1K)$ and $(1+\frac 1K)$, so it is very close to one.

\bigskip
In order to determine the permutation $\si$ of $[d-1]$ for the sequence of lines $L_1,\ldots,L_{d-1}$ we have to check, for all pairs $h<i$, the sign of
\[
\det(a_i-a_h,v_1,\ldots,v_{d-1})=\det(a_i,v_1,\ldots,v_{d-1})-\det(a_h,v_1,\ldots,v_{d-1}).
\]
Observe that we only use $\ga_i$ for special columns so it suffices to choose $\ga_{d+1},\ga_{2(d+1)},\ldots,\ga_{n(d+1)}.$
We introduce the notation $\ga_x^*=\ga_{x(d+1)}$

\begin{claim}The sequence $\ga_1^*,\ldots, \ga_n^*$ can be chosen so that  for all $h<i$
\begin{equation}\label{eq:key}
|\det(a_h,v_1,\ldots,v_{d-1})|< \eps |\det(a_i,v_1,\ldots,v_{d-1})|.
\end{equation}
\end{claim}

{\bf Proof.} The vector $v_i$ resp. $v_h$ comes from a special column $C_{x(d+1)}$ and $C_{y(d+1)}$ with $y<x$.
Then $v_i$ is in fact $v_{x(d+1)}$ and the index $i$ can be any number in $[d-1]$ except 1 because $h<i$.
Similarly $v_h$ coincides with $v_{y(d+1)}$ and $h$ can be any number in $[d-1]$ except $d-1$.

We define $\ga^*_x$ recursively, starting with $\ga^*_1=1$. Assume $\ga^*_z$ has been defined for all $z<x$, is a positive integer,
and satisfies (\ref{eq:key}).
The possible values of $|\!\det(a_i,v_1,\ldots,v_{d-1})|$ disregarding the factors $\ga^*_x$ and $(1\pm(d+1)\eps)$ (the latter is between $(1-\frac 1K)$ and $(1+\frac 1K)$) are of the form
\[
P [b_i]_{d+1-i}[(v_i)_{d+1-i}]_i
\]
for all available choices of $i$ and $v_1,\ldots, v_{d-1}$ with $v_i=v_{x(d+1)}$, $P$ also varies. This is a finite set $Z_x$ (say) of positive integers.
The possible values of $|\!\det(a_h,v_1,\ldots,v_{d-1})|$ disregarding the factor $(1\pm(d+1)\eps)$ is of the form
\[
\ga^*_yP [b_h]_{d+1-h}[(v_h)_{d+1-h}]_h
\]
for all choices of $y<x$, $h<d$ and $v_1,\ldots, v_{d-1}$ with $v_h=v_{y(d+1)}$, $P$ varies again.
This is another finite set, $V_x$ (say), of positive integers.

It is clear that there is an integer $\ga^*_x$ so that $\max V_x < \frac {\eps}{10} \min Z_x$. Bringing back the factors $1\pm(d+1)\eps$
finishes the proof. \qed

\bigskip
The claim shows that the sign of $\det(a_i-a_h,v_1,\ldots,v_{d-1})$ coincides with that of $\det(a_i,v_1,\ldots,v_{d-1})$ when $h<i$.
The sign of the last expression is the same as the sign of $(-1)^{i-1}\de_{d+1-i}$ because of (\ref{eq:det}). We can determine this sign by
choosing $\de_{d+1-i}$ any way we like. This means that when $L_{x(d+1)}$ is in position $i$ of the sequence $L_1,\ldots,L_{d-1}$,
in the corresponding permutation $\si$ of $[d-1]$, either every $h<i$ will come before $i$, or every $h<i$ will come after $i$, depending
on the choice of $\de_{d+1-i}$. This is exactly what two-sided stacked means, finishing the proof of Theorem~\ref{th:types}.\qed

\medskip
There are exactly two universal permutations for $d=3$ and six universal permutations for $d=4$, because in these cases $(d-1)!$ and $2^{d-1}-2$ coincide.
For $d=5$ the number of universal permutations is at least 14 and at most 24. The unresolved cases are the permutations 2143, 2413, 1324, 4132, 1423 and their reverses 3412, 3142, 4231, 1423, 4132, of course none of them are two-sided stacked. We do not know if any of them are universal, but we believe none of them are.

\section{Proof of Theorem~\ref{th:RI}}\label{sec:RI}

We define $a(i,j)$ by induction, and for that we require a new condition, slightly stronger than (\ref{eq:RI}) which is described next.

Consider a $k \times k$ submatrix ($k\ge 2$) $M^*$ of $M$ whose bottom right entry is $a(i,j)$. Let $P$ be the product of the elements on its main diagonal. A {\sl diagonal}, $\Delta$, of $M^*$ is a set of entries containing exactly one entry from every row and column of $M^*$. Write $\prod \Delta$ for the product of the elements in $\Delta$. We are going to require that for every diagonal $\Delta$ except the main one
\begin{equation}\label{eq:ind}
\prod \Delta < \eta P.
\end{equation}
Here we choose $\eta>0$ to be small, namely $\eta=\eps/d!$. As $\det M^*$ is the sum of all $\prod \Delta$, each taken with a well-defined sign $\pm1$, $\det M^*= P+\sum (\pm1)\prod \Delta$
where the sum is taken over all diagonals except the main one. It follows that $|\det M^* -P| \le (k!-1)\eta P$. The choice $\eta=\eps/d!<\eps/(k!-1)$ ensures that $\det M^*=(1\pm \eps)P$, exactly condition (\ref{eq:RI}).

The point in the following induction argument is that conditions (\ref{eq:ind}) and (\ref{eq:extraC}) plus the RI condition only give lower bounds on the next entry $a(i,j)$, and there are finitely many such lower bounds. This is always easy to satisfy by choosing $a(i,j)$ large enough.

Observe that (\ref{eq:ind}) only matters when $i\ge 2$. For $i=1$ we set $a(1,1)=1$ and $a(1,2)=2$. When $i=1$ condition (\ref{eq:extraC}) says that
\[
a(1,j-2)a(1,j-1)<\eps a(1,j)
\]
and only needed for $j=3,4,\ldots,m$ and is easy to satisfy by induction on $j$.

So $i\ge 2$ and we assume $a(I,J)$ have been determined and satisfy both (\ref{eq:ind}) and (\ref{eq:extraC}) for all pairs $(I,J)$ when $I<i$ and $J\in [m]$,
and when $I=i$ and $J<j$. The cases $j=1,2$ are simple but need special treatment.

Defining $a(i,1)$ is easy: the only condition is $a(i,1)>a(i-1,m)$ so $a(i,1)=a(i-1,m)+1$ will do. Defining $a(i,2)$ is similarly easy:
it has to be larger than $a(i,1)$ and condition (\ref{eq:ind}) is meaningful only for $2\times 2$ submatrices whose bottom right entry is $a(i,2)$. In this case it requires that $a(i,1)a(I,2)<\eta a(I,1)a(i,2)$ for all $I<i$. These are only lower bounds on $a(i,2)$,
and there are finitely many of them. So choosing $a(i,2)$ large enough we are done with $a(i,1)$ and $a(i,2)$.

We assume from now on that $j\ge 3$ (and $i\ge 2$). Consider a $k \times k$ ($k\ge 2$) submatrix $M^*$ of $M$
whose bottom right entry is $a(i,j)$. Write $Q$ for product of the elements
on the main diagonal of $M^*$ except $a(i,j)$. The requirement is that $\eta Q$ be larger than $\prod \Delta$ for every (but the main)
diagonal $\Delta$ of $M^*$. We consider two cases separately.

\medskip
{\bf Case 1} when $a(i,j)\notin \Delta.$ For such diagonals $\prod \Delta$ is a positive integer, independent of $a(i,j)$.
There are finitely many such integers, corresponding to every possible choice of $k\in \{2,3,\ldots,i\}$ and $M^*$ and $\Delta$, their maximum is an integer $H(i,j)$.
Moreover, let $Q(i,j)$ be the minimum of the product of the entries on the main diagonal of $M^*$ except $a(i,j)$
for every possible choice of $k$ and $M^*$. Condition (\ref{eq:ind}) requires that $H(i,j)<\eta a(i,j)Q(i,j)$.
We can clearly choose $a(i,j)$ so large that this inequality is satisfied.

\medskip
{\bf Case 2} when $a(i,j)\in \Delta.$ For these diagonals $\prod \Delta = a(i,j) \prod \Delta^*$ where $\Delta^*$ is
the corresponding diagonal of the $(k-1)\times (k-1)$ submatrix $M^{\circ}$ that you get after deleting row $i$ and
column $j$ from $M^*$.

\medskip
Assume first that $k\ge 3$. Let $a(i^{\circ},j^{\circ})$ be the bottom right entry on the main diagonal $\Delta^{\circ}$ of $M^{\circ}$.
Condition (\ref{eq:ind}) of the induction hypothesis for the pair $a(i^{\circ},j^{\circ})$ implies that
$\prod \Delta^*<\eta \prod \Delta^{\circ}$ for all diagonals of $M^{\circ}$ except the main one.
Multiplying by $a(i,j)$ we conclude that for all diagonals $\Delta$ of $M^*$ with $a(i,j) \in \Delta$ (except the main one)
$\prod \Delta <\eta a(i,j)\prod \Delta^{\circ}= \eta Q$. This shows that condition (\ref{eq:ind}) is automatically satisfied for the pair $(i,j)$ in Case 2 when $k\ge 3$.

Assume now that $k=2$. Condition (\ref{eq:ind}) says now that for all $I<i$ and $J<j$,
\[
a(i,J)a(I,j)<\eta a(I,J)a(i,j).
\]
This is again a finite set of lower bounds in $a(i,j).$ Let $H^*(i,j)$ be the maximum of these lower bonds.

 It is easy to deal with (\ref{eq:extraC}) which, in the present case $j\ge 3$, requires that
\[
a(i,j-1)a(i,j-2) < \eps a(i,j).
\]
Finally we choose an integer $a(i,j)$ larger than the maximum of the three numbers
\[
\frac {H(i,j)}{\eta Q(i,j)},\; H^*(i,j),\; \frac {a(i,j-1)a(i,j-2)}{\eps}.
\]

\section{Continuous and universal families of lines}\label{sec:cont}

In this section and the next we work with oriented lines in $\R^d$.  Let $I \subset \R$ be an open interval, for instance $I=(0,1)$ will do. Suppose that, for each $t \in I$, $L(t) = (a(t), v(t))$ is an (oriented) line in $\R^d$. We assume that $v(t)$ is a (Euclidean) unit vector, and $a(t)$ is the intersection of $L(t)$ and the hyperplane $\{x\in \R^d: v(t)\cdot x=0\}$. 

The space of oriented lines in $\R^d$ is a Grassmannian manifold with a well defined topology which defines the continuity of the map $t \to L(t)$. We say that $L(t)$ is a {\sl continuous} family of lines if $t \to L(t)$ is continuous. It is not hard to see that in this case both $a(t)$ and $v(t)$ are continuous. We will use the following fact. For a continuous family $L(t)\,(t\in I)$ of lines
\begin{equation}\label{eq:cont}
	 \lim _{s \rightarrow t} \dist( L(s), L(t))=0 
\end{equation}
where $\dist$ is the standard Euclidean distance, the infimum of the distances between any two points from the lines.

The family $L(t)$ is {\sl universal} if the type of the lines $L(t_1),\ldots,l(t_{d-1})$ is the same permutation $\si$ of $[d-1]$ for every $t_1<\ldots<t_{d-1}$, and $\si$ is the {\sl type} of such a family.

A simple example is the family of lines $L(t)=(a(t),v(t))$ with $t\in (0,1)$, say, given by $a(t)=(t,0,\ldots,0)$ and $v(t)=(1,t,\ldots,t^{d-1})$. When $0< t_1<\ldots <t_{d-1}$, suitably translated copies of the hyperplane $H=\frac 1{d-1}(L(t_1)+\ldots +L(t_{d-1}))$ contain the lines $L(t_i)$ in the order $L(t_1),\ldots,L(t_{d-1})$. This is easy to see: the translated copy of $H$ containing the line $L(t_i)$ passes through the point $a(t_i)$ and these points come in order $a(t_1),\ldots,a(t_{d-1})$ on the line $L=\{a(s):s\in \R\}.$ So this is a continuous and universal family of lines and their type is the identity. But in this example all the lines intersect the line $L$. A more generic example comes next. Some preparations are needed.

Let $0<t_1< \ldots < t_{d-1}$ be an increasing sequence of real numbers. The Vandermonde determinant
\begin{align*}
V_0 =
\begin{vmatrix}
1 & 1 & \cdots & 1 \\
t_1 & t_2 & \cdots & t_{d-1} \\
\vdots  & \vdots  & \ddots & \vdots  \\
t^{d-2}_1 & t^{d-2}_2 & \cdots & t^{d-2}_{d-1}
\end{vmatrix}
= \prod_{i=2, i>j}^{d-1}  (t_i-t_j)
\end{align*}
in positive.
Following the terminology of \cite{Heineman29} we say that $V_0$ is the {\sl principal}
Vandermondian. The generalized Vandermonde determinant is, assuming $0<b_1< \ldots < b_{d-1}$,
\begin{align*}
\begin{vmatrix}
t^{b_1}_1 & t^{b_1}_2 & \cdots & t^{b_1}_{d-1} \\
t^{b_2}_1 & t^{b_2}_2 & \cdots & t^{b_2}_{d-1} \\
\vdots  & \vdots  & \ddots & \vdots  \\
t^{b_{d-1}}_1 & t^{b_{d-1}}_2 & \cdots & t^{b_{d-1}}_{d-1}
\end{vmatrix}.
\end{align*}
This is the principal Vandermondian when $\{ b_1, \ldots, b_{d-1} \}$ coincides with $\{ 0, \ldots, d-2 \}$.
When $1 \le j\le d-1$ and $\{ b_1, \ldots, b_{d-1} \} =  \{ 0,1, \ldots, d-1 \}\setminus \{d-1-j\}$
we say that the corresponding generalized
Vandermonde determinant is the {\sl secondary} Vandermondian $V_j$.

Observe that if $j=0$ would be allowed, the definition would give back the principal Vandermondian $V_0$.

Let $E_j = \sum t_1t_2 \cdots t_j$ be the elementary symmetric function where we add
up all the products of any $j$ different variables.
Observe that $E_0 = 1$ and that $E_1$ is the sum of all variables, and $E_j>0$ since all $t_i$ are positive. \\

Theorem 1 from \cite{Heineman29} states that $\frac{V_j}{V_0} = E_j$ and therefore we have
\begin{fact}\label{fact:sv}
The secondary Vandermondians $V_1, \ldots, V_{d-1}$ are all positive.
\end{fact}

Assume now that $0<t_1 < \ldots < t_{d-1}$ and $a_1, a_2, \ldots, a_d >0$ and define the matrix
\[
A=\begin{pmatrix}
a_1 & 1 & 1 & \cdots & 1 \\
-a_2 & t_1 & t_2 & \cdots & t_{d-1} \\
\vdots & \vdots  & \vdots  & \ddots & \vdots  \\
(-1)^{d-1} a_d & t^{d-1}_1 & t^{d-1}_2 & \cdots & t^{d-1}_{d-1}
\end{pmatrix}.
\]
The determinant of $A$ when expanding by the first column is

\begin{equation}\label{eq:Vander}
\det A=a_1 V_{d-1} + a_2 V_{d-2} + \ldots + a_d V_0 > 0.
\end{equation}

We can give now the more generic example of a continuous and universal family of lines. Assume $a_i(t)>0$ is a continuous and increasing function on $t \in (0,\infty)$ for $i\in [d]$. For $t>0$ set
\begin{itemize}
	\item $v(t) = (1, t, t^2, \ldots, t^{d-1})$,
	\item $a(t) = (a_1(t), -a_2(t), \ldots, (-1)^{d-1} a_d(t) )$.
\end{itemize}

\begin{theorem}
The family of lines $L(t) = (v(t),a(t))$, $t>0$ is continuous and universal and its type is the identity.
\end{theorem}

The proof follows directly from (\ref{eq:Vander}).

\section{Uniqueness of continuous type}\label{sec:unique}

Assume $L(t)= (a(t), v(t)),\;(t \in I)$ is a continuous family of lines in $\rr.$ Here $I \subset \R$ is an open interval (possibly infinite). Recall that $v(t)$ is a unit vector for all $t\in I.$ We call the identity permutation and its reverse {\sl trivial}. The previous section contains examples of continuous and universal families of lines whose types are trivial. The target in this section is to prove the converse, namely, that the type of such a family of lines is always trivial, at least when $d\ge 5$.

\begin{theorem}\label{t:cont-unique}
Let $d\ge 5$. If $L(t)$ is a continuous and universal family of lines, then its type is trivial.
\end{theorem}

The same holds for $d=3$: the type of a continuous and universal family of lines in $\R^3$ is trivial because this is the only type available in $\R^3$. We believe that the same holds in $\R^4$ but our method does not work in that case: $\R^4$ is too small to accommodate five linearly independent vectors. 

\smallskip
Before the proof we need a few auxiliary lemmas. 
By definition, if $L(t)$ is a continuous and universal family of lines, then for any $t_1<\ldots < t_{d-1}$ the vectors $v(t_1),\ldots, v(t_{d-1})$ are linearly independent. The following claim extends this to any $d$ such vectors.

\begin{lemma}\label{l:indep} The vectors  $v(t_1), \ldots, v(t_d)$ are linearly independent when $t_1<\ldots < t_d.$
\end{lemma}

\begin{corollary}\label{cor:indep}
Let $L(t)$ be a continuous family of lines, $U$ be proper linear subspace of  $\R^d$ and $J$ some open interval $J \subset I.$ Then either $L(t)$ has trivial type, or
there exists a sub-interval $J' \subset J$ such that $v(t) \not \in U$ for all $t \in J'$.
\end{corollary}

For the next lemma we need something like the derivative of $v(t)$.
Unfortunately $v(t)$ may not be a differentiable function.

For any $t,s\in I$ let $\delta v(t,s) = v(s)-v(t)$ and define the normalized change as
$\Delta v(t,s) = \frac{\delta v(t,s)}{|| \delta v(t,s)||}$; note that $\delta v(t,s)\ne 0$ because $v(t)$ and $v(s)$ are linearly independent.
Since $\Delta v(t,s)$ lies in a compact set (namely the unit sphere), it has limit, $v'(t)$, from the right for every $t$.
That is for every $t$ there exists a sequence $s_i>t$ with limit $t$ such that the limit of $\Delta v(t,s_i)$ is $v'(t)$, and $v'(t)\ne 0.$
(There could be several different values of $v'(t)$ that could work, we just choose one of them.)

Observe that $v(t)$ and $v'(t)$ are orthogonal and that for every $t,s\in I$ with $t\ne s$ 
$\lin\{ v(t), v(s) \} = \lin \{ v(t), \Delta v(t,s) \}.$

\begin{lemma}\label{l:first-two-choices}
Let $L(t)$ be a continuous and universal family of lines in $\R^d$, $d\ge 5$.
Then either $L(t)$ has trivial type, or there exist $t,s \in I$ with $t<s$ such that the four vectors $v(t), v'(t), v(s), a(s)-a(t)$ are linearly independent.
\end{lemma}
\medskip
{\bf Proof} of Theorem~\ref{t:cont-unique} using the previous lemmas.
Let $\sigma$ be the type of $L(t)$ and assume on the contrary that $\si$ is not the identity or its reverse.
Then there exist integers $j,k$ such that $\sigma(j) < \sigma(k) < \sigma(j+1)$ or
$\sigma(j+1) < \sigma(k) < \sigma(j)$. We can assume without loss of generality that the first case occurs.

By Lemma~\ref{l:first-two-choices} we can choose
$t_j, t_k$ such that the four vectors $v(t_j), v'(t_j), v(t_k)$ and $a(t_k)-a(t_j)$
are linearly independent.
Next we choose $t_i$ one by one for $i \in [d-1] \setminus \{ j,j+1,k \}$
so that $t_1< \ldots <t_j<t_{j+2}<\ldots < t_{d-1}$
and that the $d$ vectors 	$v(t_1), \ldots, v(t_j), v'(t_j), v(t_{j+2}), \\ \ldots , v(t_{d-1}), a(t_k)-a(t_j)$
are linearly independent. For this we use Corollary~\ref{cor:indep} by defining $U$ to be the linear span of the previous (at
most $d-1$) vectors and choose $J$ to reflect the relative position of the next $t_i$.

Let $s_1 > s_2 > \ldots$ be an infinite decreasing sequence with limit $t_j$ and
$\lim_{i \to \infty}\Delta v(t_j,s_i)=v'(t_j)$. We can assume that $t_j < s_i < t_{j+2}$ for every $i.$

Set $H_i = \lin \{ v(t_1), \ldots, v(t_j), v(s_i), v(t_{j+2}), \ldots, v(t_{d-1}) \}$ which is a $(d-1)$-dimensional subspace because these vectors are linearly independent. It is clear that $H_i$ remains unchanged if in its definition $v(s_i)$ is replaced by $\Delta v(t_j,s_i).$ The orthogonal vector to $H_i$ is
\[
w_i = v(t_1) \wedge \ldots \wedge  v(t_j) \wedge \Delta v(t_j,s_i) \wedge
 v(t_{j+2}) \wedge \ldots \wedge v(t_{d-1}). 
\]
Let $u_i = \frac{w_i}{||w_i||}$ be the unit vector orthogonal to $H_i$.

The limit of $H_i$ is $H = \lin \{ v(t_1), \ldots, v(t_j), v'(t_j), v(t_{j+2}), \ldots, v(t_{d-1}) \}$ and the limit of $u_i$ is $u$, a unit vector orthogonal to $H$.

The order of the lines $L(t_1), \ldots, L(t_j), L(s_i), L(t_{j+2}), \ldots, L(t_{d-1})$
is the same as the order of the real numbers
$$u_i \cdot a(t_1),\ldots,  u_i \cdot a(t_j),
 u_i \cdot a(s_i), u_i \cdot a(t_{j+2}),\ldots,  u_i \cdot a(t_{d-1}).$$
 Since $L(s_i)$ plays the role of the $(j+1)$st line
 $$
  u_i \cdot a(t_j) <  u_i \cdot a(t_k) <  u_i \cdot a(s_i)
 $$
The distance of the two lines $L(t_j)$ and $L(s_i)$ is at least
the distance of any two hyperplanes through the two lines, particularly the ones parallel to $H_i$, which is $u_i \cdot (a(s_i)-a(t_j)).$ It follows from (\ref{eq:cont}) that, as $i$ tends to infinity, the limit of $ u_i \cdot a(t_j)$ and $ u_i \cdot a(s_i)$ is the same, namely $u \cdot a(t_j)$.
Therefore the limit of the previous inequality is
 $$
u \cdot a(t_j) \le  u \cdot a(t_k) \le  u \cdot a(t_j)
$$
and we have $u \cdot (a(t_k)-a(t_j)) =0$. Consequently the unit vector $u$ is orthogonal to every one of the $d$ vectors $$v(t_1), \ldots, v(t_j), v'(t_j), v(t_{j+2}), \ldots , v(t_{d-1}), a(t_k)-a(t_j)$$ that are linearly independent.
A contradiction. \qed

\bigskip
{\bf Proof} of Lemma~\ref{l:indep}. Assume that the vectors $v(t_1),\ldots,v(t_d)$ are not linearly independent. Their linear span is then a $(d-1)$-dimensional subspace. Let $z$ be its normal vector. The numbers $h_i=z\cdot a(t_i),\; (i\in [d])$ come in the order $h_{i_1},...,h_{i_d}$, and $\pi=(i_1,...,i_d)$ is a permutation of $[d].$ We can assume that $i_1<i_d$ (by replacing $z$ by $-z$ if necessary). The permutation $\si_j$ of $[d-1]$ comes from $\pi$ by deleting the entry $j\in [d]$. Universality implies that, for every $j,$ $\si_j$ is the {\sl same ordering} of $d-1$ linearly ordered elements. Assume $d=i_h$, then in $\si_j$ the largest element is in position $h-1$ if $j \in \{i_1,...,i_{h-1}\}$ and in position $h$ if $j\in \{i_{h+1},...,i_d\}$ implying that $\{i_{h+1},...,i_d\}=\emptyset.$ It follows that $h=d.$ The same argument works for the $k$th largest element of $[d]$ (by backward induction on $k$): it has to be in position $d-k+1$ in $\pi$. \qed

\medskip
{\bf Proof} of Corollary~\ref{cor:indep}. Let $S \subset I$ be the set of all $t$, such that $v(t) \in U$.
Observe that $|S| \le d-1$ as otherwise $U$ coincides with $\rr$ because any $d$ direction vectors are linearly independent by Lemma~\ref{l:indep}.
We can choose $J\subset I$ to be any open interval avoiding $S.$\qed

\medskip
{\bf Proof} of Lemma~\ref{l:first-two-choices}.
Assume on the contrary that the four vectors $v(t), v'(t), v(s), a(s)-a(t)$ are linearly dependent for every pair $t<s$. Set

\[
W(t,s)=\lin \{v(t), v'(t), v(s), a(s)-a(t)\}.
\]
The vectors $v(t), v(s), a(s)-a(t)$ are linearly independent because any two lines of the family are skew. Consequently $\dim W(t,s)=3$ and $W(t,s) = \lin \{ v(t), v(s), a(s)-a(t)\}.$

Given $t_0<t_1<t_2$ set
\[
W = \lin \{ v(t_0), v'(t_0), v(t_1), a(t_1)-a(t_0), v(t_2), a(t_2)-a(t_0) \}.
\]
Observe that $W$ is the span of the two $3$-dimensional subspaces $W(t_0,t_1)$ and $W(t_0,t_2)$
that intersect in an at least  $2$-dimensional subspace because both contain $\lin \{ v(t_0), v'(t_0) \}$.
Thus the dimension of $W$ is at most $4$.

Moreover $W$ contains $W(t_1,t_2)=\lin \{ v(t_1), v(t_2), a(t_2)-a(t_1)\}$ as these three vectors all belong to $W$. Then $v'(t_1)\in W$ since $v'(t_1)\in W(t_1,t_2).$

Assume next that there exist $t_0<t_1$ such that the four vectors $v(t_0), v'(t_0), v(t_1), v'(t_1)$ are linearly independent. Then there exists $t_2>t_1$ such that $v(t_0), v'(t_0), v(t_1), v'(t_1), v(t_2)$ are linearly independent because otherwise $v(t_2)$ lies in $\lin \{v(t_0), v'(t_0), v(t_1), v'(t_1), v(t_2) \}$ for every $t_2>t_1$ contradicting Lemma~\ref{l:indep}. It follows that the linear span of $v(t_0), v'(t_0), v(t_1), v'(t_1), v(t_2)$ is $5$-dimensional. But it is a subspace of $W$ and $\dim W\le 4,$ a contradiction showing that the four vectors $v(t_0), v'(t_0), v(t_1), v'(t_1)$ are indeed linearly dependent for every $t_0< t_1.$ Note that this is the point where the condition $d>4$ is used.

So any two of the $2$-dimensional subspaces $\lin \{ v(t), v'(t) \}$ intersect in an at least  $1$-dimensional subspace.
There are two possibilities that can happen.

\smallskip
{\bf Case 1}. All subspaces $\lin \{ v(t), v'(t) \}$ are contained in some $3$ dimensional space.
This is impossible as there exist $d\ge 4$ direction vectors that are independent.\\

{\bf Case 2}. There exists some vector $u$ contained in $\lin \{v(t), v'(t) \}$ for every $t$.
In this case choose an increasing sequence $t_0<\ldots<t_n$ such that
$u$ together with any $d-1$ different direction vectors $v(t_i)$ are linearly independent.
This goes by induction on $n$ and for $n\le d-1$ a simple application of Corollary~\ref{cor:indep} works. When $n \ge d$ and we have defined $t_0,\ldots,t_{n-1}$, then $t_n$ is found by repeated, actually $n-1 \choose d-2$-fold, applications of Corollary~\ref{cor:indep}. 

\smallskip
For every $t_i$, $i\in [n]$, define $\ga_i \in \R$ by
$$
a(t_i)-a(t_0) = \ga_i u + \alpha_i v(t_i) + \beta_i v(t_0)
$$
where $\alpha_i, \beta_i$ are some real numbers. The existence of $\ga_i,\al_i,\be_i$ follows from the fact that  $W(t_0,t_i)=\lin \{v(t_0),v(t_i),a(t_i)-a(t_0)\}$ and $u\in \lin \{v(t_0),v'(t_0)\} \subset W(t_0, t_i)$. Note that $\ga_i$ remains the same even if $a(t_i)\in L(t_i)$ is replaced by another point $a(t_i)+\de v(t_i)\in L.$

Observe that
\begin{eqnarray*}
a(t_i)-a(t_j) &=& a(t_i)-a(t_0) - (a(t_j)-a(t_0) )\\
  &=& (\ga_i-\ga_j) u + \alpha_i v(t_i) -\alpha_j v(t_j) + (\beta_i-\beta_j) v(t_0).
\end{eqnarray*}

Assume $n>d^2$. A classic result of Erd\H os and Szekeres~\cite{ESz35} shows then that there exists a subsequence of the $t_i$ of length $(d-2)$ such that the corresponding $\ga_i$ form an increasing (or decreasing) sequence. For simplicity assume that this subsequence is $t_1, \ldots, t_{d-2}$ and the $\ga_i$ are increasing. Let $w = v(t_0) \wedge v(t_1) \wedge  \ldots \wedge v(t_{d-2}) $ be the orthogonal vector to all lines $L(t_i)$. Then the order of the lines $L(t_i)$ depends on the order of the numbers $w \cdot a(t_i)$
and therefore on the signs of $w \cdot (a(t_i)-a(t_j)) = (\ga_i-\ga_j) w \cdot u$.

So the order of all the lines but $L(t_0)$ is increasing, meaning that in $\si$ the last $d-2$ elements come in increasing order. But now by symmetry, choose $t_n$ as the anchor of the $\ga_i$s and we get that the first $d-2$ elements in $\si$ come in increasing order. Therefore $\si$ is identity permutation.\qed

\bigskip
{\bf Acknowledgements.}  Research of IB was partially supported by Hungarian National Research Grants (no. 131529, 131696, and 133819), and research of GK by the Israel Science Foundation (grant no. 1612/17).

\bigskip
\noindent
Imre B\'ar\'any \\
R\'enyi Institute of Mathematics,\\
13-15 Re\'altanoda Street, Budapest, 1053 Hungary\\
{\tt barany.imre@renyi.hu} and\\
Department of Mathematics\\
University College London\\
Gower Street, London, WC1E 6BT, UK

\medskip
\noindent
Gil Kalai\\
Einstein Institute of Mathematics\\
Hebrew University,
Jerusalem 91904, Israel,\\
{\tt kalai@math.huji.ac.il} and\\
Efi Arazy School of Computer Science,
IDC, Herzliya, Israel

\medskip
\noindent
Attila P\'or\\
Department of Mathematics\\
Western Kentucky University\\
1906 College Heights Blvd. \#11078\\
Bowling Green, KY 42101, USA\\
{\tt  attila.por@wku.edu}\\

\end{document}